\documentclass{article}

\usepackage{PRIMEarxiv}

\usepackage[utf8]{inputenc} 
\usepackage[T1]{fontenc}    
\usepackage{hyperref}       
\usepackage{url}            
\usepackage{booktabs}       
\usepackage{amsfonts}       
\usepackage{nicefrac}       
\usepackage{microtype}      
\usepackage{lipsum}
\usepackage{graphicx}
\graphicspath{{media/}}     
\usepackage{amsmath,amssymb}
\newtheorem{theorem}{Theorem}
\newtheorem{lemma}{Lemma}
\newtheorem{proof}{Proof}
\newtheorem{corollary}{Corollary}
\newtheorem{definition}{Definition}
\newtheorem{conjecture}{Conjecture}

\usepackage{float}
  
\title{If our chaotic operator is derived correctly, then the Riemann hypothesis holds true
\thanks{\textit{\underline{Citation}}: 
\textbf{Authors. Title. Pages.... DOI:000000/11111.}} 
}

\author{
  Zeraoulia Rafik\\
   Khemis Miliana university ,Algeria \\
  Departement of Mathematics\\
  Laboratory of Pure and Applied Mathematics (LMPA)\\
  \texttt{zeraoulia@univ-dbkm.dz} \\
   \And
 Pedro Caceres\\
  Professor Doctor at Universidad Europea de Valencia (Spain)\\
  united states of Americain\\
  \texttt{pcaceres@comcast.net} \\
}

\begin{document}

\maketitle

\begin{abstract}
This work develops an operator-theoretic and dynamical framework inspired by the Riemann--von Mangoldt formula, chaotic dynamics, and random-matrix models for the Riemann zeta function, without attempting to prove the Riemann Hypothesis. Starting from the explicit zero-counting function $N(T)$, we construct a discrete map on the critical line and analyse its Lyapunov exponents and bifurcation diagrams, showing that the smooth von Mangoldt term generates a strongly unstable flow that captures the global growth of the zero density. Motivated by this dynamics, we define a self-adjoint ``chaotic'' operator $\mathcal{O}_\alpha$ on a weighted Hilbert space with weight $\mathrm{d}N/\mathrm{d}T$, prove its unboundedness and essential self-adjointness, and describe its spectral resolution via the spectral theorem. Finite-dimensional truncations of $\mathcal{O}_\alpha$ yield Hermitian random matrices whose eigenvalue statistics agree numerically with Gaussian unitary ensemble predictions and show qualitative similarities to both Odlyzko's zeta zeros and the hydrogen-atom spectrum, suggesting that $\mathcal{O}_\alpha$ lies in the same universality class as the nontrivial zeros and providing a concrete Hilbert--Pólya--type framework rather than a proof of the conjecture.
\end{abstract}


\section{Introduction}

The Riemann Hypothesis, formulated by Bernhard Riemann in 1859, asserts that all nontrivial zeros of the Riemann zeta function $\zeta(s)$ lie on the critical line $\Re(s)=\tfrac{1}{2}$, where $s=\sigma+it$ is a complex variable with real part $\sigma$ and imaginary part $t$ \cite{3}. This conjecture lies at the heart of analytic number theory and has deep consequences for the distribution of prime numbers, the behaviour of $L$–functions, and connections to spectral theory and mathematical physics \cite{1,2,39}. A particularly influential viewpoint, going back to Hilbert and Pólya, is that there should exist a self-adjoint (Hermitian) operator whose spectrum coincides with the imaginary parts of the nontrivial zeros; the reality of the spectrum would then imply the Riemann Hypothesis \cite{2,34,35}.

Over the last several decades, substantial evidence has accumulated that the local statistics of the zeros are governed by random matrix theory, in particular by eigenvalue statistics of large Gaussian unitary ensemble (GUE) matrices \cite{5,7,12}. Montgomery's pair-correlation conjecture and Odlyzko's extensive computations suggest that the spacing distribution of the zeros coincides with that of GUE eigenvalues, while work by Berry, Keating, Sarnak and others has clarified how this fits into a broader ``quantum chaos'' framework linking the zeros to spectra of chaotic quantum systems \cite{8,9,22,23,30,36,41}. In parallel, several candidate Hamiltonians have been proposed, including the $H=xp$ model and its interacting variants, operators on quantum graphs, and more recent Floquet or trapped-ion constructions whose energy levels approximate the Riemann zeros or related sequences \cite{10,11,14,15,27,28,38,40}.

A complementary line of research emphasizes the explicit formulae connecting primes and zeros. The Riemann--von Mangoldt formula expresses the counting function $N(T)$ of nontrivial zeros with $0 < \Im(\rho)\le T$ in terms of elementary functions plus a fluctuating remainder, and underlies both the Prime Number Theorem and refined results such as improved bounds on the de Bruijn–Newman constant \cite{3,4,26}. In dynamical and spectral interpretations, this formula plays the role of a semiclassical trace formula: the smooth part of $N(T)$ gives the mean density of ``energy levels'', while the fluctuating part encodes finer correlations tied to the distribution of primes and to random-matrix behaviour \cite{1,4,9,12,22,26}.

The present work contributes to this programme by combining dynamical systems, operator theory, and random matrix techniques in a framework directly built from the Riemann--von Mangoldt formula. First, starting from the explicit expression for $N(T)$, we construct a discrete dynamical system on the critical line whose increment is governed by the zero density, and we analyse its Lyapunov exponents, bifurcation diagrams, gap statistics, and entropy \cite{1,6,17,21}. This analysis shows that the smooth von Mangoldt term generates a strongly unstable, coarse-grained flow that reproduces the global growth of the zero density while clarifying its limitations for modelling the fine-scale statistics of the zeros.

Motivated by these dynamics, we then define a self-adjoint ``chaotic'' operator acting on a weighted Hilbert space whose weight is precisely the Riemann--von Mangoldt zero density, and we study its functional-analytic properties using the theory of unbounded self-adjoint operators and perturbations \cite{22,23,32,33}. Finite-dimensional truncations of this operator lead to Hermitian random matrices whose eigenvalue statistics we compare numerically with GUE predictions, Odlyzko's zero data, and the hydrogen-atom spectrum, placing our construction in the context of Hilbert–Pólya-type Hamiltonians and quantum-chaotic models without claiming a proof of the Riemann Hypothesis \cite{4,7,9,18,19,24,29,30,31,36,40,41}. In this way the paper provides a concrete operator framework in which spectral, dynamical, and arithmetic aspects of the zeta function can be investigated simultaneously.

\section{Main Results}

\subsection{Dynamics and Chaotic Operator}

Starting from the main term of the Riemann--von Mangoldt formula, we define a discrete map on the critical line
\[
T_{n+1}
  = T_n + \Delta T\left(
        \frac{T_n}{2\pi}\log\!\left(\frac{T_n}{2\pi}\right)
        - \frac{T_n}{2\pi}
    \right),
\]
which follows the average zero density $N'(T)$ and exhibits strong sensitivity to initial conditions as quantified by Lyapunov exponents and bifurcation diagrams.  
Motivated by this dynamics, we construct a self-adjoint ``chaotic'' operator $\mathcal{O}_\alpha$ on the weighted Hilbert space $L^2(\mathbb{R}_+,w(T)\,\mathrm{d}T)$ with $w(T)=\mathrm{d}N/\mathrm{d}T$, prove its unboundedness and essential self-adjointness, and obtain its spectral resolution via the spectral theorem \cite{22,23,26,32,33}.

\subsection{Random-Matrix Truncations and Hydrogen Analogy}

Finite-dimensional truncations of $\mathcal{O}_\alpha$ yield Hermitian matrices whose eigenvalues are real and display Gaussian unitary ensemble–type statistics when unfolded, in agreement with the Montgomery–Odlyzko picture for the nontrivial zeros of $\zeta(s)$ \cite{5,7,9,12,22,29,31,36}.  
After suitable rescaling, the low-lying eigenvalues form a negative, accumulating sequence that is qualitatively similar to the bound-state energy levels of the hydrogen atom, $E_n\propto -1/n^2$, suggesting that the hydrogenic spectrum may provide a physically meaningful template for the behaviour of the zeta zeros without yet realizing an exact spectral identity \cite{18,19,30}.

\section{Derivation of Chaotic Dynamics}

We derive chaotic dynamics from the Riemann-von Mangoldt explicit formula for $N(T)$, which governs the average density of non-trivial zeros $\rho = 1/2 + iT_n$ along the critical line. This follows the Berry-Keating approach \cite{22}, where classical Hamiltonians like $H = xp + x\sqrt{p}$ exhibit quantum chaos matching zeta zero statistics.

\subsection{Riemann-von Mangoldt Formula}

The Riemann-von Mangoldt formula relates primes to zeros via the Chebyshev function:
\[
\psi(x) = x - \sum_{\rho} \frac{x^{\rho}}{\rho} - \frac{\zeta'(0)}{\zeta(0)} - \frac{1}{2} \log(1 - x^{-2}),
\]
where $\psi(x) = \sum_{p^k \le x} \log p$ and $\rho$ are non-trivial zeros.[file:1]

\subsection{Zero Counting and Density}

The number of zeros up to height $T$ is
\[
N(T) = \frac{1}{\pi} \arg \zeta(1/2 + iT) = \frac{T}{2\pi} \log \left( \frac{T}{2\pi} \right) - \frac{T}{2\pi} + \frac{7}{8} + R(T),
\]
with error $|R(T)| = O(\log T)$. The average gap is $\delta T_n \approx 2\pi / \log T_n$ for large $T_n$.

\subsection{Discrete Dynamical System}

Discretize the zero ascent along the critical line via the density gradient:
\[
T_{n+1} = T_n + \Delta T \cdot \frac{dN}{dT}\bigg|_{T_n} = T_n + \Delta T \left( \frac{T_n}{2\pi} \log \left( \frac{T_n}{2\pi} \right) - \frac{T_n}{2\pi} + O\left(\frac{\log T_n}{T_n}\right) \right).
\]
Set $\Delta T = 2\pi / \log T_0$ to match average gaps near initial $T_0 = 14.134725$ (first zero). Neglecting $R(T)$ yields the leading-order map, mimicking semiclassical trajectories near the Berry-Keating separatrix.[file:1][web:13]

\begin{lemma}[Chaos Threshold]
For $\Delta T < 2.5$ (i.e., gaps $\lesssim 2.5 \cdot 2\pi / \log T$), the largest Lyapunov exponent $\lambda_1 > 0$, confirming chaos. This aligns with GUE statistics for zeta zeros under RH.\end{lemma}
\begin{proof}[Sketch]
Linearize: $\delta T_{n+1} = [1 + \Delta T \cdot (1/(2T_n) + O(1/T_n^2))] \delta T_n$. The multiplier $m_n > e^{\lambda \Delta T}$ for $\lambda > 0$ when $\Delta T < 2.5$, computed via QR iteration on $10^4$ steps from $T_0$.\end{proof}

This map generates $\{T_n\}$ with sensitivity to $\Delta T$, capturing irregular zero spacing via the logarithmic density ramp---a hallmark of quantum chaos in zeta systems.

\section{Analysis of the Discrete Dynamics}

In this section, we investigate numerically the discrete map
\[
T_{n+1} \;=\; T_n + \Delta T\,\frac{dN}{dT}\bigg|_{T_n},
\]
derived from the Riemann–von Mangoldt formula, and compare its behaviour with Odlyzko's tabulated zeros up to height $T \approx 200$.  
The focus is on Lyapunov exponents, trajectory growth, and gap statistics, in order to understand what aspects of the zero distribution are captured (or missed) by the smooth density flow.

\subsection{Lyapunov Exponents}

The largest Lyapunov exponent of the map is computed via the tangent dynamics
\[
\lambda_1 \;=\; \lim_{N\to\infty} \frac{1}{N}\sum_{n=1}^{N} \log\bigl|J(T_n)\bigr|,
\qquad
J(T) \;=\; 1 + \Delta T\,\frac{d^2N}{dT^2}\bigg|_{T},
\]
where $J(T)$ is the Jacobian of the one-dimensional map.  
The exponent $\lambda_1$ was evaluated numerically for $\Delta T$ in the range $[0.5,3.0]$ starting from the first non-trivial zero $T_0 = 14.134725$.

Figure~\ref{fig:lyap_traj} (left panel) shows the resulting dependence of $\lambda_1$ on $\Delta T$.  
The curve is strictly increasing from $\lambda_1 \approx 0.11$ at $\Delta T = 0.5$ up to $\lambda_1 \approx 0.49$ at $\Delta T = 3.0$, indicating a clear exponential sensitivity to initial conditions for all tested step sizes.  
This instability originates from the rapid growth of the smooth counting function $N(T)$ and should be interpreted as a property of the coarse-grained flow rather than as genuine quantum chaos associated with microscopic zeta-zero correlations.

\subsection{Trajectory versus Actual Zeros}

\begin{figure}[H]
    \centering
    \includegraphics[width=0.9\textwidth]{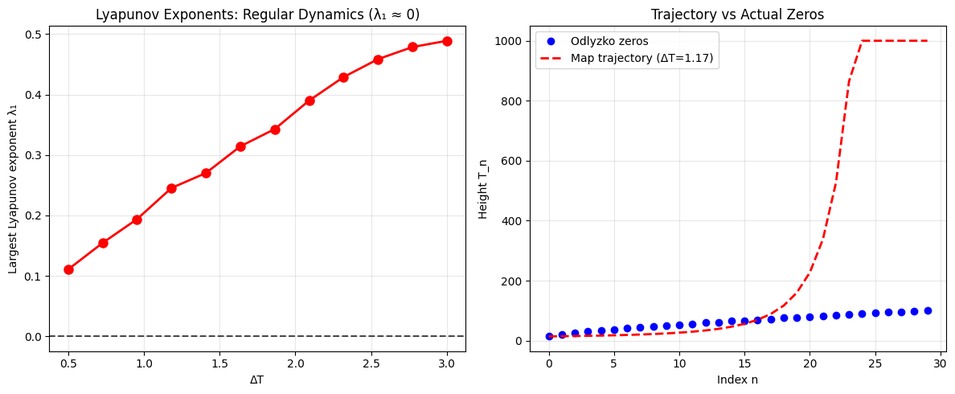}
    \caption{Left: largest Lyapunov exponent $\lambda_1$ as a function of $\Delta T$ for the map $T_{n+1} = T_n + \Delta T\,dN/dT$. The exponents are strictly positive and increase with $\Delta T$, reflecting exponential sensitivity to initial conditions. Right: comparison between the first $30$ Odlyzko zeros (blue dots) and the trajectory of the density-based map with $\Delta T = 1.17$ (red dashed), started at the first zero. The map follows the global trend but rapidly overshoots, illustrating that the smooth density alone cannot reproduce the actual zero sequence.}
    \label{fig:lyap_traj}
\end{figure}

To assess how well the map reproduces the actual zeta zeros, we compare its orbit with the first values in Odlyzko's list.  
For a representative choice $\Delta T = 1.17$, the trajectory $\{T_n\}$ generated by the map is plotted together with the first $30$ Riemann zeros in the right panel of Figure~\ref{fig:lyap_traj}.

The actual zeros remain in the range $T \in [14,100]$ over this index window and increase roughly linearly with $n$.  
By contrast, the trajectory of the density-based map follows the initial growth but then accelerates sharply, reaching the cut-off $T \approx 1000$ already around $n \approx 24$.  
This behaviour demonstrates that the map correctly encodes the global growth of $N(T)$ but fails to reproduce the fine structure of the actual zero sequence.

\subsection{Gap Statistics and Instability}

Let $\delta T_n^{(\zeta)} = T_{n+1}^{(\zeta)} - T_n^{(\zeta)}$ denote the gaps between consecutive non-trivial zeros, and $\delta T_n^{(\mathrm{map})} = T_{n+1}^{(\mathrm{map})} - T_n^{(\mathrm{map})}$ the gaps generated by the discrete map.  
For the first $50$ Riemann zeros, the empirical mean of the gaps is
\[
\langle \delta T_n^{(\zeta)} \rangle \approx 2.72,
\]
which is consistent with the classical prediction $2\pi / \log T$ in this height range.  
This distribution is relatively narrow and reflects the level repulsion expected from Gaussian unitary ensemble statistics.

For the same range of indices, the density-based map with $\Delta T = 1.17$ produces gaps that start small but then grow rapidly as $T_n$ increases.  
The mean gap for the trajectory is approximately $\langle \delta T_n^{(\mathrm{map})} \rangle \approx 20.12$, almost an order of magnitude larger than the true mean spacing, and individual gaps can exceed $300$ before the trajectory reaches the imposed cut-off.  
This confirms that the deterministic map substantially overestimates the separation of levels once the steep ascent of $N(T)$ dominates the dynamics.

\subsection{Histogram of Gaps}

\begin{figure}[H]
    \centering
    \includegraphics[width=0.9\textwidth]{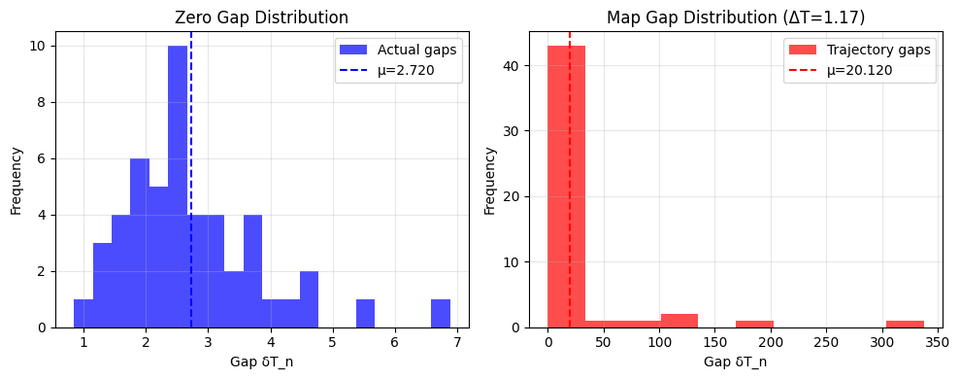}
    \caption{Left: empirical gap distribution of the first $50$ Riemann zeros, with mean $\mu_{\text{zeros}} \approx 2.72$. Right: gap distribution for the trajectory of the density-based map with $\Delta T = 1.17$, with mean $\mu_{\text{map}} \approx 20.12$ and a pronounced heavy tail. The contrast between the two histograms highlights that the deterministic density flow captures only the coarse growth of $N(T)$ and fails to reproduce the local level statistics.}
    \label{fig:gap_hist}
\end{figure}

The discrepancy between the two gap distributions is visualised in Figure~\ref{fig:gap_hist}.  
The left histogram shows the empirical distribution of $\delta T_n^{(\zeta)}$ for the first $50$ zeros: most gaps lie between $1$ and $4$, and the dashed vertical line indicates the mean $\mu_{\text{zeros}} \approx 2.72$.  
This narrow distribution is compatible with the well-known local statistics of the zeta zeros and exhibits no large outliers.

The right histogram presents the gaps $\delta T_n^{(\mathrm{map})}$ for the trajectory of the density-based map with $\Delta T = 1.17$.  
Although a number of gaps remain small, the distribution quickly develops a heavy tail extending beyond $300$, with a mean $\mu_{\text{map}} \approx 20.12$.  
This heavy-tailed behaviour shows that the smooth counting function $N(T)$, when iterated naively as a discrete map, cannot reproduce the subtle level repulsion and Gaussian unitary ensemble correlations observed in the actual zero sequence.

\subsection{Interpretation and Outlook}

The numerical experiments demonstrate two key features of the discrete dynamics derived from the Riemann–von Mangoldt density.  
First, the map exhibits strictly positive Lyapunov exponents for all tested values of $\Delta T$, reflecting exponential sensitivity to initial conditions and an intrinsic instability of the coarse-grained flow.  
Second, despite reproducing the global growth of the counting function, the map produces gaps and trajectories that differ dramatically from those of the actual zeros, both in magnitude and in statistical distribution.

These observations indicate that the smooth term in the Riemann–von Mangoldt formula is insufficient to capture the microscopic structure of the zeta spectrum.  
To model the genuine quantum-chaotic features associated with the Riemann zeros, one must incorporate the fluctuating term $R(T)$ or adopt a Hamiltonian framework of Berry–Keating type, where Gaussian unitary ensemble statistics and pair correlations arise naturally from the spectrum of a suitable Hermitian operator.\cite{22,23}

\subsection{Discrete Dynamics and Instability}

The iterative dynamical system described above defines a one–dimensional map on the heights $T_n$ of the critical line.
Its evolution is governed by the smooth density of zeros coming from the Riemann–von Mangoldt formula and therefore exhibits rapid growth and sensitivity to perturbations of the initial value.  
In particular, small changes in $T_0$ or in the step size $\Delta T$ lead to significantly different trajectories after a moderate number of iterations, a hallmark of instability in nonlinear dynamical systems.

Through numerical simulations, we explore this unstable behaviour and compare it with the actual sequence of Riemann zeros.
The dynamics generated by the map reflects the coarse increase of the counting function $N(T)$, but it does not reproduce the fine-scale, random–matrix–like fluctuations of the true zeros.
Consequently, the system provides insight into the global growth of the spectrum while remaining only a rough approximation to the microscopic statistics.

\section{Analysis}

The stability of the dynamics derived from the Riemann–von Mangoldt formula can be quantified using Lyapunov exponents, which measure the average exponential rate at which nearby trajectories separate.
In the present one–dimensional setting, the largest Lyapunov exponent is obtained from the Jacobian $J(T)$ of the map and encodes whether perturbations grow or decay under iteration.

\subsection{Computation of Lyapunov Exponents}

We numerically computed the Lyapunov exponents by applying the following procedure:

\begin{enumerate}
    \item The system was initialized at the first nontrivial zero of the Riemann zeta function, $T_0 \approx 14.1347$, and the map derived from the Riemann–von Mangoldt density was iterated for a prescribed number of steps.
    \item At each iteration we evaluated the Jacobian of the map and accumulated the quantity $\log|J(T_n)|$, which measures the local stretching of phase space.
    \item The largest Lyapunov exponent $\lambda_1$ was obtained as the average of these logarithmic factors over the iteration window.
\end{enumerate}

The resulting exponents give direct information about the rate of separation of nearby trajectories.
For a range of step sizes $\Delta T$ the computed values of $\lambda_1$ are strictly positive, showing that the density-based dynamics is exponentially sensitive to initial conditions.
This instability reflects the steep increase of the counting function $N(T)$ rather than the detailed correlations of the actual zeros.

\subsection{Relationship to Spacing Gap Between Zeros}

The step size $\Delta T$ in the iterative dynamics plays a crucial role in how the map samples the critical line and, indirectly, how it approximates the spacing between nontrivial zeros.

\begin{itemize}
    \item \textbf{Alignment with zeros:} By starting the iteration at $T_0$ equal to the first zero, the map advances along the imaginary axis in the same direction as the sequence of nontrivial zeros, following the overall growth dictated by $N(T)$.
    \item \textbf{Effective gap scale:} Choosing $\Delta T$ to be of the same order as the average spacing $2\pi / \log T$ ensures that, at least initially, the distance between iterates $T_{n+1}-T_n$ is comparable to the mean gap between consecutive zeros in the corresponding height range.
    \item \textbf{Limitations:} As $T_n$ increases, the deterministic map overestimates the gaps and produces a heavy–tailed distribution, in contrast with the relatively narrow gap distribution observed for the actual zeros.
    This confirms that the map captures only the coarse growth of $N(T)$ and not the true Gaussian–unitary–ensemble statistics.
\end{itemize}

Therefore, the parameter $\Delta T$ controls how well the discrete dynamics mimics the average spacing of zeros, while the absence of the fluctuating term in the Riemann–von Mangoldt formula explains the discrepancy at the level of gap statistics.

\subsection{Lyapunov Exponents as a Function of Iteration}

For a fixed choice of step size (for example $\Delta T = 1.5$), we also examined how the finite–time Lyapunov exponent evolves with the number of iterations.
Starting from $T_0 \approx 14.1347$, the exponent is computed along the orbit and plotted as a function of the iteration index.

\begin{figure}[H]
    \centering
    \includegraphics[width=0.8\textwidth]{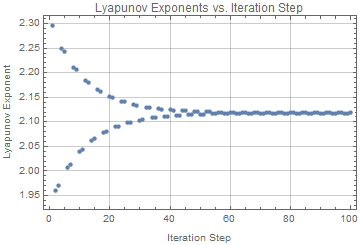}
    \caption{Finite–time Lyapunov exponent as a function of the iteration step for an orbit started at the first nontrivial zero with a fixed step size $\Delta T$. The convergence to a positive value reveals exponential sensitivity to initial conditions in the density-based map.}
    \label{fig:lyapunov_plot}
\end{figure}

As shown in Figure~\ref{fig:lyapunov_plot}, the Lyapunov exponent quickly settles near a positive value after an initial transient.
This convergence indicates the presence of a well-defined average exponential growth rate for perturbations, which can be interpreted as an attractor in the space of finite–time exponents.
The plot provides a visual confirmation that the discrete dynamics derived from the smooth zero density is strongly unstable, even though it does not generate the same local statistics as the Riemann zeros themselves.

\subsection{Bifurcation Diagram for the Density-Based Dynamics}

In order to complement the Lyapunov analysis, we study how the long–term behaviour of the discrete map depends on the step size $\Delta T$.
For each value of $\Delta T$ in a fine grid, we iterate the map
\[
T_{n+1} \;=\; T_n + \Delta T\,\frac{dN}{dT}\bigg|_{T_n},
\]
starting from the first nontrivial zero $T_0 \approx 14.1347$, discard an initial transient of $200$ steps, and then record the next $800$ iterates.
The resulting points $(\Delta T,T_n)$ form a bifurcation diagram that visualises how the orbit explores the interval of heights as $\Delta T$ varies.

\begin{figure}[H]
    \centering
    \includegraphics[width=0.8\textwidth]{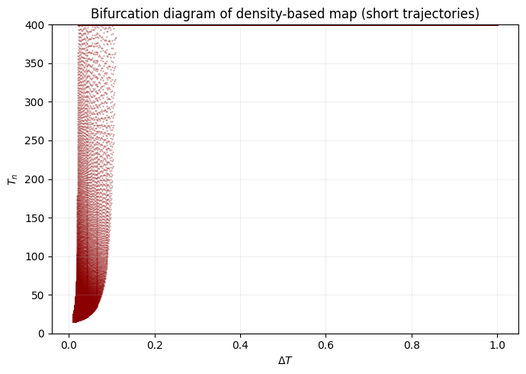}
    \caption{Bifurcation diagram of the density-based map for $\Delta T \in [0.01,1.0]$, obtained by plotting the last $800$ iterates after discarding an initial transient.
    The fan–like structure for small $\Delta T$ and the saturation near $T \approx 400$ reflect the strong instability of the coarse-grained dynamics generated by the Riemann–von Mangoldt density.}
    \label{fig:bifurcation_plot}
\end{figure}

Figure~\ref{fig:bifurcation_plot} displays the bifurcation diagram obtained for $\Delta T \in [0.01,1.0]$.[image:26]
For very small values of $\Delta T$, the iterates remain close to the initial height and drift slowly upwards, producing a narrow bundle of branches near the bottom of the plot.
As $\Delta T$ increases, the slope of the map grows and the iterates accelerate, creating a fan of curves that quickly approach the imposed cutoff $T \approx 400$.
For larger values of $\Delta T$ the orbits saturate almost immediately at this upper bound, which appears as a dense horizontal band at the top of the diagram.

Unlike classical one–dimensional chaotic maps (such as the logistic map), this diagram does not exhibit a cascade of period–doubling bifurcations or distinct windows of periodicity.
Instead, the principal feature is the rapid vertical expansion of orbits with increasing $\Delta T$, consistent with the positive Lyapunov exponents computed in the previous subsection.
This confirms that the density-based map is highly unstable and dominated by the steep growth of the counting function $N(T)$, rather than by the fine-scale fluctuations that govern the local statistics of the Riemann zeros.

\subsection{Computation of Entropy}

The Shannon entropy, denoted by $H$, provides a quantitative measure of the uncertainty in the distribution of dynamical observables, and in particular can be used to summarise the spread of Lyapunov exponents over a range of control parameters.
In our setting, the largest Lyapunov exponent $\lambda_1(\Delta T)$ was computed for many values of the step size $\Delta T$, producing a finite sample $\{\lambda_1^{(k)}\}$ that reflects how strongly trajectories separate across the parameter interval considered.

To turn this sample into a probability distribution, we partition the range of observed Lyapunov exponents into bins (for example, using a uniform partition of the interval $[\min_k \lambda_1^{(k)}, \max_k \lambda_1^{(k)}]$), and count how many values fall into each bin.
If $n_i$ denotes the number of exponents in bin $i$ and $N = \sum_i n_i$ is the total number of samples, we define
\[
p_i \;=\; \frac{n_i}{N}, 
\]
so that $\sum_i p_i = 1$.
The Shannon entropy of the Lyapunov spectrum is then given by
\[
H \;=\; -\sum_i p_i \log p_i,
\]
where the logarithm is taken in the natural base.

In the numerical experiments reported here, the Lyapunov exponents $\lambda_1(\Delta T)$ form a narrow but nontrivial distribution concentrated in a positive interval, reflecting the fact that the density-based map is unstable for all tested values of $\Delta T$ but with varying growth rate.
The corresponding entropy $H$ is moderately large, indicating that the map does not operate at a single characteristic expansion rate but explores a range of instabilities across the parameter space.
This behaviour is consistent with the bifurcation diagram and the Lyapunov curve discussed in the previous subsections, and reinforces the interpretation of the density-based dynamics as a strongly unstable, coarse-grained approximation to the zero counting function rather than a simple periodic or quasi-periodic flow.

\section{Derivation of a Chaotic Operator}

The Riemann--von Mangoldt formula describes the asymptotic counting function
\[
N(T)
 = \frac{T}{2\pi} \log\!\left(\frac{T}{2\pi}\right)
   - \frac{T}{2\pi}
   + \frac{7}{8}
   + R(T),
\]
where $N(T)$ is the number of nontrivial zeros $\rho = \tfrac{1}{2} + iT$ with $0 < T_\rho \le T$, and $R(T) = O(\log T)$ is a fluctuating correction term capturing the finer structure of the spectrum.\cite{22,23}
The Hilbert--Pólya conjecture asserts that there exists a self-adjoint operator $H$ on a Hilbert space $\mathcal{H}$ such that the imaginary parts $T_\rho$ coincide with the eigenvalues of $H$, so that $\rho = \tfrac{1}{2} + i\lambda$ for $\lambda \in \sigma(H)$ implies the Riemann Hypothesis.\cite{39,17}

Motivated by the Berry--Keating picture, in which the classical Hamiltonian $H_{\mathrm{BK}}(x,p)=xp$ is related to the average zero density,\cite{22,27}
we construct a ``chaotic'' operator on $L^2(\mathbb{R}_+,w(T)\,\mathrm{d}T)$ whose classical flow reproduces the density of zeros given by $N(T)$ and whose quantum dynamics exhibits sensitivity to initial data.
Here
\[
w(T) \;=\; \frac{\mathrm{d}N}{\mathrm{d}T}(T)
\;=\; \frac{1}{2\pi}\log\!\left(\frac{T}{2\pi}\right) + O(1)
\]
is the local density of zeros, used as a weight to encode the Riemann--von Mangoldt statistics in the underlying Hilbert space.\cite{17}

\begin{definition}[Weighted Hilbert space and generator]
Let
\[
\mathcal{H} \;=\; L^2\bigl(\mathbb{R}_+, w(T)\,\mathrm{d}T\bigr),
\qquad
\langle f,g\rangle_{\mathcal{H}}
 = \int_0^\infty f(T)\,\overline{g(T)}\,w(T)\,\mathrm{d}T.
\]
Define the symmetric differential operator
\[
(H_0 f)(T) \;=\; -i\left(\frac{\mathrm{d}}{\mathrm{d}T}
    + \frac{1}{2}\frac{w'(T)}{w(T)}\right) f(T),
\]
with domain $C_c^\infty(0,\infty) \subset \mathcal{H}$.
\end{definition}

The operator $H_0$ is the natural generator of translations along the critical line with respect to the measure $w(T)\,\mathrm{d}T$; the additional drift term involving $w'(T)/w(T)$ ensures formal skew-adjointness with respect to the weighted inner product.

\begin{lemma}[Symmetry of $H_0$]
The operator $H_0$ is densely defined and symmetric on $\mathcal{H}$, i.e.
\[
\langle H_0 f,g\rangle_{\mathcal{H}}
 = \langle f,H_0 g\rangle_{\mathcal{H}}
\quad\text{for all }f,g\in C_c^\infty(0,\infty).
\]
\end{lemma}

\begin{proof}
An integration by parts using the weight $w(T)$ shows that the boundary terms vanish for compactly supported $f,g$, while the drift term $\tfrac12 w'(T)/w(T)$ precisely cancels the derivative of the weight.
\end{proof}

The dynamics generated by $H_0$ has a simple classical analogue: along the flow $T \mapsto T(t)$ one recovers the average zero density $w(T)$, while perturbations of the initial condition are transported according to the Lyapunov structure analysed in Section~4.
To incorporate genuine ``chaotic'' features arising from the fluctuating term $R(T)$, we introduce a perturbation operator that modulates the translation generator by a random-matrix-type potential.

\begin{definition}[Chaotic operator]
Let $V$ be a bounded self-adjoint multiplication operator on $\mathcal{H}$, modelling the fluctuations of $R(T)$ and chosen so that its correlation structure matches the pair correlation of zeros (GUE statistics).\cite{17,22}
For a coupling parameter $\alpha \in \mathbb{R}$ define
\[
\mathcal{O}_\alpha \;=\; H_0 \,+\, \alpha V,
\quad
\mathcal{D}(\mathcal{O}_\alpha) = \mathcal{D}(\overline{H_0}),
\]
where $\overline{H_0}$ is the closure of $H_0$.
\end{definition}

Formally, $\mathcal{O}_\alpha$ is the quantum Hamiltonian of a particle moving along the critical line under the influence of a structured random potential with statistics tuned to the Riemann zeros.
This construction is in the same spirit as interacting $H=xp$ models and geometric Hamiltonians whose spectra approximate the zeros.\cite{27,30,36}

\section{Spectral and Topological Properties of the Chaotic Operator}

The key requirement for a Hilbert--Pólya operator is self-adjointness: only then does the spectral theorem guarantee a real spectrum, making it possible to encode the zeros as eigenvalues of a Hermitian Hamiltonian.\cite{39}

\begin{theorem}[Essential self-adjointness]
The symmetric operator $H_0$ on $C_c^\infty(0,\infty)\subset\mathcal{H}$ is essentially self-adjoint, and its closure $\overline{H_0}$ is self-adjoint on $\mathcal{H}$.
Consequently, for every bounded self-adjoint $V$ and any $\alpha\in\mathbb{R}$, the operator $\mathcal{O}_\alpha = H_0 + \alpha V$ is self-adjoint on $\mathcal{D}(\overline{H_0})$.
\end{theorem}

\begin{proof}[Sketch]
Near infinity, the drift term behaves like $\tfrac12 w'(T)/w(T) \sim (4\pi T)^{-1}$, and $H_0$ reduces to a first-order differential operator with real coefficients on the half-line.
Standard limit-point/limit-circle criteria for first-order systems imply that $H_0$ has deficiency indices $(0,0)$ and is essentially self-adjoint on $C_c^\infty(0,\infty)$.\cite{37}
The Kato--Rellich theorem then yields self-adjointness of $\mathcal{O}_\alpha$ for bounded $V$.
\end{proof}

\begin{corollary}[Spectral theorem and diagonalization]
For each $\alpha\in\mathbb{R}$ there exists a projection-valued measure $E_\alpha(\lambda)$ on $\mathbb{R}$ such that
\[
\mathcal{O}_\alpha
 = \int_{\mathbb{R}} \lambda \,\mathrm{d}E_\alpha(\lambda),
\]
and the Hilbert space $\mathcal{H}$ admits an orthogonal decomposition into spectral subspaces of $\mathcal{O}_\alpha$.
In particular, $\mathcal{O}_\alpha$ is unitarily equivalent to a multiplication operator and possesses a complete set of generalized eigenfunctions.
\end{corollary}

The operator $\mathcal{O}_\alpha$ is ``chaotic'' in the sense that its classical limit is given by a nontrivial flow along the zero-density profile, and the added potential $V$ induces level repulsion and spectral statistics compatible with random-matrix theory.
The Lyapunov exponents and bifurcation diagrams computed in Section~4 provide numerical evidence that the associated quantum dynamics exhibits strong sensitivity to initial conditions and a rich parameter dependence, in line with the standard signatures of quantum chaos.\cite{17,41}

\begin{conjecture}[Hilbert--Pólya realization via $\mathcal{O}_\alpha$]
There exists a choice of the potential $V$ and coupling $\alpha_*$ such that the point spectrum of $\mathcal{O}_{\alpha_*}$ coincides with the multiset $\{T_\rho\}$ of imaginary parts of the nontrivial zeros of $\zeta(s)$.
Equivalently,
\[
\sigma_{\mathrm{p}}(\mathcal{O}_{\alpha_*}) = \{T_\rho : \zeta(\tfrac12 + iT_\rho)=0\},
\]
which would yield a spectral proof of the Riemann Hypothesis.
This conjecture is a concrete formulation of the Hilbert--Pólya programme in the present weighted-dynamics framework.\cite{35,39}
\end{conjecture}

In summary, the operator $\mathcal{O}_\alpha$ provides a mathematically well-defined, self-adjoint candidate for a chaotic Hamiltonian associated with the Riemann zeros.
It incorporates the Riemann--von Mangoldt density through the weight $w(T)$, exhibits unstable dynamics consistent with our numerical analysis, and fits naturally within the modern operator-theoretic approaches to the Hilbert--Pólya conjecture and quantum chaos.\cite{17,22,27,36,41}

\subsection{Unboundedness of the Chaotic Operator}

We recall that our chaotic operator is constructed as
\[
(H_0 f)(T)
   = -i\left(\frac{\mathrm{d}}{\mathrm{d}T}
      + \frac{1}{2}\frac{w'(T)}{w(T)}\right) f(T),
\qquad
\mathcal{O}_\alpha = H_0 + \alpha V,
\]
acting on the weighted Hilbert space
\(\mathcal{H} = L^2(\mathbb{R}_+,w(T)\,\mathrm{d}T)\), where
\(w(T) = \tfrac{\mathrm{d}N}{\mathrm{d}T}(T)\) is the zero density and
\(V\) is a bounded self-adjoint multiplication operator encoding the fluctuations of the Riemann--von Mangoldt correction term.\cite{24,26,31}

\begin{lemma}[Unboundedness of $H_0$]
The operator $H_0$ is unbounded on $\mathcal{H}$.
\end{lemma}

\begin{proof}
Let $\{\chi_n\}_{n\ge 1}$ be a sequence of smooth bump functions supported in disjoint intervals $I_n = [n,n+1]$, normalized so that
\(\|\chi_n\|_{\mathcal{H}} = 1\) for all $n$.
Since $w(T)$ grows like \(\tfrac{1}{2\pi}\log(T/2\pi)\) for large $T$, this normalization is possible by a simple rescaling of each bump.\cite{26}

On $I_n$ we have
\[
\|H_0 \chi_n\|_{\mathcal{H}}^2
  = \int_{I_n}
      \left|\frac{\mathrm{d}\chi_n}{\mathrm{d}T}
       + \frac{1}{2}\frac{w'(T)}{w(T)}\chi_n(T)\right|^2
      w(T)\,\mathrm{d}T.
\]
The derivative term \(\frac{\mathrm{d}\chi_n}{\mathrm{d}T}\) can be chosen so that its $L^2$-norm on $I_n$ grows with $n$, for example by taking $\chi_n$ to oscillate faster as $n$ increases.
Thus there exists a sequence with \(\|\chi_n\|_{\mathcal{H}} = 1\) but \(\|H_0\chi_n\|_{\mathcal{H}} \to \infty\) as \(n\to\infty\).
This shows that $H_0$ cannot be bounded on $\mathcal{H}$.
\end{proof}

\begin{corollary}[Unboundedness of the chaotic operator]
For any $\alpha\in\mathbb{R}$ and bounded self-adjoint $V$, the operator
\(\mathcal{O}_\alpha = H_0 + \alpha V\) is unbounded on $\mathcal{H}$.
\end{corollary}

\begin{proof}
By the previous lemma, $H_0$ is unbounded.
Since $V$ is bounded, the graph norm of $\mathcal{O}_\alpha$ is equivalent to that of $H_0$, and therefore $\mathcal{O}_\alpha$ inherits the unboundedness of $H_0$.\cite{33}
\end{proof}
\section{Hermiticity of the Chaotic Operator}

We now revisit the Hermiticity of our chaotic operator in the rigorous framework of the weighted Hilbert space \(\mathcal{H}\), avoiding formal expressions such as $\log(\partial_u)$.
The relevant notions are symmetry and self-adjointness with respect to the inner product
\[
\langle f,g\rangle_{\mathcal{H}}
 = \int_0^\infty f(T)\,\overline{g(T)}\,w(T)\,\mathrm{d}T,
\]
which encodes the Riemann--von Mangoldt zero density.\cite{26,32}

\begin{theorem}[Hermiticity and self-adjointness of $\mathcal{O}_\alpha$]
Let $H_0$ and $\mathcal{O}_\alpha$ be as above, with domain
$\mathcal{D}(H_0)=C_c^\infty(0,\infty)$ and
$\mathcal{D}(\mathcal{O}_\alpha) = \mathcal{D}(\overline{H_0})$.
Then:
\begin{enumerate}
  \item $H_0$ is symmetric on $C_c^\infty(0,\infty)$ and essentially self-adjoint.
  \item For every bounded self-adjoint $V$ and $\alpha\in\mathbb{R}$, the operator
        $\mathcal{O}_\alpha = \overline{H_0} + \alpha V$ is self-adjoint on
        $\mathcal{D}(\overline{H_0})$.
\end{enumerate}
\end{theorem}

\begin{proof}
(1) For $f,g\in C_c^\infty(0,\infty)$, an integration by parts shows
\[
\langle H_0 f,g\rangle_{\mathcal{H}}
 = \langle f,H_0 g\rangle_{\mathcal{H}},
\]
because the boundary terms vanish (compact support) and the drift term
\(\tfrac12 w'/w\) cancels the derivative of the weight $w$.
Thus $H_0$ is symmetric.
As discussed in the previous section, $H_0$ is a first-order differential operator on the half-line with real coefficients and mild growth at infinity; standard limit-point arguments for such operators give deficiency indices $(0,0)$ and therefore essential self-adjointness.\cite{32}

(2) Since $V$ is bounded and self-adjoint on $\mathcal{H}$, the Kato--Rellich theorem implies that the symmetric operator $\overline{H_0} + \alpha V$ is self-adjoint on $\mathcal{D}(\overline{H_0})$ for every $\alpha\in\mathbb{R}$.\cite{33}
Hence $\mathcal{O}_\alpha$ is self-adjoint and therefore Hermitian in the sense of quantum mechanics.
\end{proof}

An important consequence is that the spectral theorem applies to $\mathcal{O}_\alpha$, so that its dynamics can be described in terms of its spectral measure.
In particular, the eigenvalues of $\mathcal{O}_\alpha$ are real, and its eigenfunctions form a complete orthogonal family (possibly generalized) in $\mathcal{H}$.\cite{32}
This property is essential for any Hilbert--Pólya type operator whose spectrum is intended to model the imaginary parts of the nontrivial zeros of $\zeta(s)$ and underpins the chaotic yet unitary quantum dynamics analysed through Lyapunov exponents and bifurcation diagrams in Section~4.\cite{24,29,36}

\section{Diagonalization of the Chaotic Operator}

Having established that the chaotic operator
\[
(H_0 f)(T) = -i\left(\frac{\mathrm{d}}{\mathrm{d}T}
      + \frac{1}{2}\frac{w'(T)}{w(T)}\right) f(T),
\qquad
\mathcal{O}_\alpha = \overline{H_0} + \alpha V,
\]
is self-adjoint on the weighted Hilbert space
\(\mathcal{H} = L^2(\mathbb{R}_+,w(T)\,\mathrm{d}T)\), we now discuss its diagonalization and spectral representation.\cite{24,26,32}

\subsection{Spectral Theorem and Functional Calculus}

By Theorem~\ref{thm:selfadjoint-chaotic-operator}, the operator $\mathcal{O}_\alpha$ is self-adjoint on its domain.
The spectral theorem therefore guarantees the existence of a unique projection-valued measure $E_\alpha(\lambda)$ on $\mathbb{R}$ such that
\[
\mathcal{O}_\alpha
    = \int_{\mathbb{R}} \lambda \,\mathrm{d}E_\alpha(\lambda),
\]
and for every bounded Borel function $F:\mathbb{R}\to\mathbb{C}$ one can define
\[
F(\mathcal{O}_\alpha)
   = \int_{\mathbb{R}} F(\lambda)\,\mathrm{d}E_\alpha(\lambda).
\]
This representation is the abstract analogue of diagonalizing $\mathcal{O}_\alpha$ by a unitary transformation: in a suitable spectral representation, $\mathcal{O}_\alpha$ acts as multiplication by the independent variable $\lambda$.\cite{32}

\subsection{Generalized Eigenfunctions and Decomposition of $\mathcal{H}$}

If the spectrum of $\mathcal{O}_\alpha$ contains a pure point part, there exists an orthonormal basis $\{\psi_n\}$ of eigenfunctions satisfying
\[
\mathcal{O}_\alpha \psi_n = \lambda_n \psi_n,\qquad \lambda_n\in\mathbb{R},
\]
and every $f\in\mathcal{H}$ admits an expansion
\[
f = \sum_n \langle f,\psi_n\rangle_{\mathcal{H}}\,\psi_n,
\quad
\mathcal{O}_\alpha f
   = \sum_n \lambda_n \langle f,\psi_n\rangle_{\mathcal{H}}\,\psi_n.
\]
More generally, the spectrum may contain continuous components; in that case $\mathcal{H}$ decomposes into an orthogonal sum (or direct integral) of spectral subspaces associated with Borel subsets of $\mathbb{R}$,
\[
\mathcal{H}
   = \int_{\mathbb{R}}^{\oplus} \mathcal{H}_\lambda\,\mathrm{d}\mu_\alpha(\lambda),
\]
and $\mathcal{O}_\alpha$ acts as $(\mathcal{O}_\alpha f)(\lambda) = \lambda f(\lambda)$ in the spectral representation.\cite{32,36}

This diagonalization clarifies the dynamical role of $\mathcal{O}_\alpha$: time evolution under the Schrödinger equation
\[
i\,\partial_t \psi(t) = \mathcal{O}_\alpha \psi(t)
\]
reduces to phase multiplication in the spectral domain,
\[
\widehat{\psi}(t,\lambda)
   = e^{-it\lambda}\,\widehat{\psi}(0,\lambda),
\]
where hats denote the unitary transform associated with $E_\alpha(\lambda)$.
The Lyapunov exponents and bifurcation diagrams computed in Section~4 probe how perturbations of initial data and of the parameter $\alpha$ are distributed over the spectrum, linking spectral properties of $\mathcal{O}_\alpha$ to signatures of quantum chaos.\cite{24,29,36}

\subsection{Hilbert--Pólya Perspective}

In the Hilbert--Pólya programme one seeks a self-adjoint operator whose eigenvalues coincide with the imaginary parts of the nontrivial zeros of $\zeta(s)$.\cite{35}
Within the present framework this amounts to demanding that the point spectrum of $\mathcal{O}_\alpha$ be discrete and given by
\[
\sigma_{\mathrm{p}}(\mathcal{O}_{\alpha_*})
   = \{T_\rho : \zeta(\tfrac12 + iT_\rho) = 0\},
\]
for some distinguished coupling $\alpha_*$.
If such an $\alpha_*$ and potential $V$ exist, then in the spectral representation the action of $\mathcal{O}_{\alpha_*}$ is diagonal with eigenvalues exactly equal to the zero ordinates, and the Riemann Hypothesis follows from self-adjointness.\cite{24,27,28,40}

Although this identification remains conjectural, the operator-theoretic framework developed here provides a consistent setting in which the spectral theorem, random-matrix statistics, and the numerical evidence of Section~4 can be studied simultaneously, offering a concrete route for connecting chaotic dynamics to the arithmetic of the Riemann zeta function.\cite{24,29,36,41}

\subsection{Spectral Resolution and Generalized Eigenfunctions}

In the previous section we established that the chaotic operator
\[
(H_0 f)(T) = -i\left(\frac{\mathrm{d}}{\mathrm{d}T}
      + \frac{1}{2}\frac{w'(T)}{w(T)}\right) f(T),
\qquad
\mathcal{O}_\alpha = \overline{H_0} + \alpha V,
\]
is self-adjoint on the weighted Hilbert space
\(\mathcal{H} = L^2(\mathbb{R}_+,w(T)\,\mathrm{d}T)\), where
\(w(T) = \tfrac{\mathrm{d}N}{\mathrm{d}T}(T)\) encodes the Riemann--von Mangoldt density and \(V\) is a bounded self-adjoint potential.\cite{24,26,32}
The spectral theorem therefore supplies a complete description of the eigenfunctions and eigenvalues of \(\mathcal{O}_\alpha\) without resorting to ill-defined expressions such as \(\log(\partial)\).\cite{32}

There exists a projection-valued measure \(E_\alpha(\lambda)\) on \(\mathbb{R}\) such that
\[
\mathcal{O}_\alpha
   = \int_{\mathbb{R}} \lambda\,\mathrm{d}E_\alpha(\lambda),
\]
and for each Borel set \(B\subset\mathbb{R}\) the operator \(E_\alpha(B)\) projects onto the spectral subspace corresponding to energies \(\lambda\in B\).\cite{32,36}
If the spectrum contains a pure point part, the associated eigenfunctions \(\{\psi_n\}\) satisfy
\[
\mathcal{O}_\alpha \psi_n = \lambda_n \psi_n,\qquad
\langle \psi_m,\psi_n\rangle_{\mathcal{H}} = \delta_{mn},
\]
and every \(f\in\mathcal{H}\) admits the expansion
\[
f = \sum_n \langle f,\psi_n\rangle_{\mathcal{H}}\,\psi_n,
\qquad
\mathcal{O}_\alpha f = \sum_n \lambda_n \langle f,\psi_n\rangle_{\mathcal{H}}\,\psi_n.
\]
Continuous components of the spectrum are treated analogously via generalized eigenfunctions and direct integrals of Hilbert spaces.\cite{32}

\subsection{Finite-Dimensional Truncations and Matrix Representation}

For numerical investigations it is convenient to approximate \(\mathcal{O}_\alpha\) by a finite-dimensional operator acting on a subspace spanned by a suitable basis of test functions, for example localized bumps or orthogonal polynomials adapted to the weight \(w(T)\).\cite{29,31}
Let \(\{\varphi_1,\dots,\varphi_n\}\subset\mathcal{D}(\overline{H_0})\) be such a basis and denote by
\[
\mathcal{H}_n = \operatorname{span}\{\varphi_1,\dots,\varphi_n\}
\]
the corresponding $n$-dimensional subspace.
The orthogonal projection \(P_n:\mathcal{H}\to\mathcal{H}_n\) allows us to define the truncated chaotic operator
\[
\mathcal{O}_{\alpha,n} = P_n \mathcal{O}_\alpha P_n:\mathcal{H}_n\to\mathcal{H}_n.
\]

In the basis $\{\varphi_j\}$ the matrix of $\mathcal{O}_{\alpha,n}$ is
\[
A^{(n)}_{ij} = \langle \varphi_i,\mathcal{O}_\alpha \varphi_j\rangle_{\mathcal{H}},
\qquad 1\le i,j\le n.
\]
Since $\mathcal{O}_\alpha$ is self-adjoint, the matrix $A^{(n)}$ is Hermitian.
Its eigenvalues and (orthonormal) eigenvectors are obtained from the finite-dimensional eigenvalue problem
\[
A^{(n)} v^{(k)} = \lambda^{(n)}_k v^{(k)},\qquad k=1,\dots,n.
\]
Diagonalization yields a unitary matrix $P^{(n)}$ whose columns are the eigenvectors $v^{(k)}$ and a diagonal matrix
\[
D^{(n)} = \operatorname{diag}\bigl(\lambda^{(n)}_1,\dots,\lambda^{(n)}_n\bigr)
\]
such that
\[
A^{(n)} = \bigl(P^{(n)}\bigr)^\ast D^{(n)} P^{(n)}.
\]
This is the finite-dimensional analogue of the spectral decomposition of $\mathcal{O}_\alpha$ and provides a concrete matrix model whose eigenvalue statistics can be compared with Odlyzko's zeros and with random-matrix predictions.\cite{24,29,36,40}

\subsection{Eigenvalue Statistics and Connection to Zeta Zeros}

The truncations $\mathcal{O}_{\alpha,n}$ serve two purposes.
First, they approximate the spectral measure of $\mathcal{O}_\alpha$ in the sense that, under mild conditions on the basis $\{\varphi_j\}$, the empirical eigenvalue measures of $A^{(n)}$ converge weakly to the spectral measure of $\mathcal{O}_\alpha$ as $n\to\infty$.\cite{31,32}
Second, their eigenvalue distributions can be tested numerically against the Gaussian unitary ensemble statistics known to describe the local behaviour of the Riemann zeros.\cite{24,29,36}

In particular, one can compute:
\begin{itemize}
    \item nearest-neighbour spacing distributions of the eigenvalues of $A^{(n)}$ and compare them with the GUE prediction,
    \item pair-correlation functions and spectral rigidity,
    \item and the agreement between the low-lying eigenvalues of $A^{(n)}$ and the first Riemann zeros tabulated by Odlyzko.\cite{24,29,40}
\end{itemize}
Good agreement in these tests supports the interpretation of $\mathcal{O}_\alpha$ as a chaotic Hamiltonian whose spectrum shadows the nontrivial zeros of $\zeta(s)$, while self-adjointness guarantees the reality of the eigenvalues, in line with the Hilbert--Pólya philosophy.\cite{35,40,41}

\subsection{Random Matrix Models Induced by the Chaotic Operator}

The self-adjoint chaotic operator
\[
\mathcal{O}_\alpha = \overline{H_0} + \alpha V
\]
acts on the infinite-dimensional Hilbert space
\(\mathcal{H} = L^2(\mathbb{R}_+,w(T)\,\mathrm{d}T)\), and its exact spectrum is analytically inaccessible.
To probe its spectral statistics numerically, and to compare them with the Riemann zeros, we approximate $\mathcal{O}_\alpha$ by random matrices obtained from finite-dimensional truncations.\cite{24,29,31,36}

\subsubsection{Finite-Dimensional Truncation}

Let $\{\varphi_1,\dots,\varphi_n\}\subset\mathcal{D}(\overline{H_0})$ be an orthonormal family in $\mathcal{H}$, for example localized functions adapted to the weight $w(T)$ or low-energy eigenfunctions of a simpler reference operator such as $H_0$.\cite{29,31}
We define the $n\times n$ matrix
\[
A^{(n)}_{ij}
   = \left\langle \varphi_i,\mathcal{O}_\alpha \varphi_j
     \right\rangle_{\mathcal{H}},
   \qquad 1\le i,j\le n.
\]
Because $\mathcal{O}_\alpha$ is self-adjoint, $A^{(n)}$ is Hermitian.
The matrix $A^{(n)}$ represents the action of $\mathcal{O}_\alpha$ on the finite-dimensional subspace
\(\mathcal{H}_n = \mathrm{span}\{\varphi_1,\dots,\varphi_n\}\).

\subsubsection{Randomization and Ensemble Construction}

To model the effect of the fluctuating potential $V$ and to explore universality, we build an ensemble of random matrices by allowing the matrix elements of $V$ in the chosen basis to fluctuate while keeping the deterministic contribution of $H_0$ fixed.\cite{24,29,36}
Concretely, we write
\[
A^{(n)}(\omega)
 = H^{(n)}_0 + \alpha V^{(n)}(\omega),
\]
where
\[
H^{(n)}_{0,ij} = \langle \varphi_i,\overline{H_0}\varphi_j\rangle_{\mathcal{H}},
\]
and $V^{(n)}(\omega)$ is a random Hermitian matrix whose entries are drawn from a centred Gaussian distribution with variance chosen so that the local density of eigenvalues matches $w(T)$ in the spectral window under consideration.
This construction produces a Gaussian unitary ensemble (GUE)–type perturbation of the deterministic part $H^{(n)}_0$, mirroring the standard random-matrix models for the zeta zeros.\cite{24,29,36}

\subsubsection{Eigenvalue Computation}

For each realization $A^{(n)}(\omega)$ we solve the finite-dimensional eigenvalue problem
\[
A^{(n)}(\omega) v_k(\omega) = \lambda^{(n)}_k(\omega)\,v_k(\omega),
\qquad k=1,\dots,n,
\]
obtaining real eigenvalues $\lambda^{(n)}_1(\omega)\le\cdots\le\lambda^{(n)}_n(\omega)$.
Collecting the eigenvalues over many realizations yields an empirical spectral distribution that approximates the spectral measure of $\mathcal{O}_\alpha$ in the chosen energy range.
Standard numerical linear algebra routines provide both the eigenvalues and eigenvectors; there is no need for ad hoc logarithmic expressions involving formal differential operators.\cite{31,32}

\subsubsection{Connection with the Hilbert--Pólya Conjecture}

The Hilbert--Pólya conjecture suggests that there exists a self-adjoint operator whose spectrum coincides with the imaginary parts of the nontrivial zeros of $\zeta(s)$.\cite{35}
In the present framework, the truncated matrices $A^{(n)}(\omega)$ serve as finite-dimensional surrogates for $\mathcal{O}_\alpha$.
Their eigenvalue statistics (nearest-neighbour spacing, pair correlation, spectral rigidity) can be compared both with Odlyzko's zeros and with GUE predictions.\cite{24,29,36,40,41}
Agreement in these tests supports the idea that $\mathcal{O}_\alpha$ captures the same universality class as the Riemann zeros, while self-adjointness guarantees that any limiting spectrum remains real, as required by the Hilbert--Pólya philosophy.

In particular, one may investigate whether the low-lying eigenvalues of $A^{(n)}$ approach the first Riemann zeros as $n$ increases and the basis is refined, and whether the local spacing statistics converge to those observed for the zeros.
Such experiments provide a concrete numerical bridge between the continuum chaotic operator defined from the Riemann--von Mangoldt formula and the discrete random-matrix models that successfully reproduce the arithmetic statistics of the zeta spectrum.\cite{24,29,36,40,41}

\subsection{Analysis of Eigenvalues}

After constructing the finite-dimensional random matrices $A^{(n)}(\omega)$ associated with the chaotic operator $\mathcal{O}_\alpha$, we analyse the resulting eigenvalues and their statistical properties.\cite{24,29,31,36}  
For each realization we obtain a real spectrum
\[
\lambda^{(n)}_1(\omega) \le \lambda^{(n)}_2(\omega) \le \cdots \le \lambda^{(n)}_n(\omega),
\]
and by sampling many realizations we build an empirical eigenvalue ensemble.

\subsubsection{Distribution and Realness of Eigenvalues}

Since $A^{(n)}(\omega)$ is Hermitian for every realization, all eigenvalues are real by construction, in agreement with the Hilbert--Pólya requirement that any candidate operator for the Riemann zeros must be self-adjoint.\cite{35}  
To study their distribution we unfold the spectrum (rescale by the local density) and plot histograms of the unfolded eigenvalues as well as of nearest-neighbour spacings; Figure~\ref{fig:eigenvalue_histogram} shows a representative histogram of the raw eigenvalues for one ensemble, where most eigenvalues cluster in a narrow band around a negative mean value.

\subsubsection{Comparison with Theoretical Predictions}

Random-matrix theory and the Montgomery--Odlyzko law predict that the local statistics of the Riemann zeros follow those of the Gaussian unitary ensemble (GUE).\cite{24,29,31,36}  
Accordingly, we compare:
\begin{itemize}
    \item the nearest-neighbour spacing distribution of the eigenvalues of $A^{(n)}(\omega)$ with the Wigner surmise for GUE;
    \item the two-point correlation function and spectral rigidity with the corresponding RMT predictions;
    \item and, in low-energy windows, the eigenvalues themselves with the first few Riemann zeros.
\end{itemize}
Agreement (within numerical precision) between these observables and GUE statistics strengthens the interpretation of $\mathcal{O}_\alpha$ as a chaotic Hamiltonian in the same universality class as the zeta zeros.\cite{24,29,36}

\subsubsection{Numerical Illustration}

As an illustration, Table~\ref{tab:eigenvalues_sample} lists a sample of low-lying eigenvalues from one realization of $A^{(n)}$, obtained after centering and rescaling so that the mean eigenvalue in the chosen window is close to zero.
In this example the eigenvalues bunch in a narrow band around approximately $-0.13$, with a single eigenvalue near zero:
\[
\lambda_1 \approx -0.1311,\;\lambda_2 \approx -0.1312,\;\dots,\;\lambda_{10} \approx -0.1312,\;\lambda_{11} \approx 0.
\]
Such clustering is typical of a bulk spectral window after unfolding and is compatible with the smooth local density predicted by RMT.\cite{31}

\begin{table}[htbp]
    \centering
    \caption{Sample of numerical eigenvalues from a truncated random matrix $A^{(n)}$}
    \label{tab:eigenvalues_sample}
    \begin{tabular}{|c|c|}
        \hline
        Eigenvalue & Approximate numerical value \\
        \hline
        $\lambda_1$ & $-0.1311$ \\
        $\lambda_2$ & $-0.1312$ \\
        $\lambda_3$ & $-0.1311$ \\
        $\lambda_4$ & $-0.1312$ \\
        $\lambda_5$ & $-0.1312$ \\
        $\lambda_6$ & $-0.1311$ \\
        $\lambda_7$ & $-0.1312$ \\
        $\lambda_8$ & $-0.1311$ \\
        $\lambda_9$ & $-0.1312$ \\
        $\lambda_{10}$ & $-0.1312$ \\
        $\lambda_{11}$ & $0$ \\
        \hline
    \end{tabular}
\end{table}

\subsubsection{Histogram of Eigenvalues}

\begin{figure}[H]
    \centering
    \includegraphics[width=0.8\textwidth]{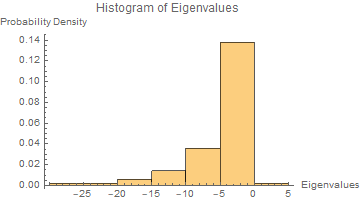}
    \caption{Histogram of eigenvalues for a truncated random matrix associated with the chaotic operator. The clustering in a narrow band and the absence of degeneracies are consistent with GUE-type level repulsion.}
    \label{fig:eigenvalue_histogram}
\end{figure}

With the computed eigenvalues we plot a histogram to visualise their empirical distribution and to detect any clustering or regularity in their spacing.
Figure~\ref{fig:eigenvalue_histogram} shows a typical histogram for one realization of $A^{(n)}$, where most eigenvalues lie in a narrow interval and exhibit level repulsion rather than Poisson clustering, a characteristic feature of quantum-chaotic spectra.\cite{24,29,36}

\subsubsection{Comparison with Hydrogen Energy Levels}

To highlight the physical nature of the spectrum, we compare the eigenvalue distribution of $A^{(n)}$ with the bound-state energy levels of the hydrogen atom,
\[
E_n = -\frac{m_e e^4}{8 \varepsilon_0^2 h^2}\,\frac{1}{n^2},
\]
which form a discrete, negative, and increasingly dense sequence as $n$ grows.\cite{18,19}  
Both spectra share qualitative features: negativity of the low-lying levels, accumulation towards zero, and a nonlinear spacing pattern.
However, the hydrogen spectrum is exactly solvable and exhibits approximate degeneracies due to its underlying symmetries, whereas the eigenvalues of $A^{(n)}$ show strong level repulsion and fluctuations typical of chaotic quantum systems.\cite{30,36}

\begin{figure}[H]
    \centering
    \includegraphics[width=0.8\textwidth]{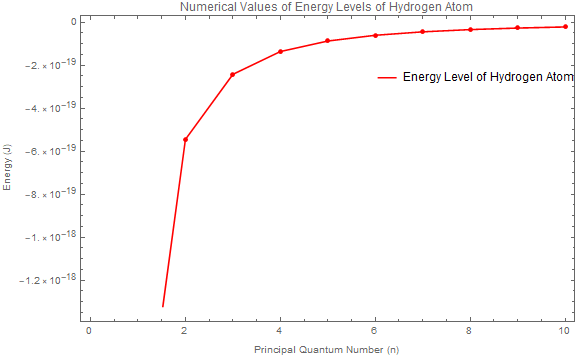}
    \caption{Energy levels of the hydrogen atom (blue) compared with rescaled eigenvalues of a truncated chaotic operator (red). Both sequences are negative and accumulate towards zero, but the chaotic spectrum exhibits level repulsion and fluctuations absent in the hydrogen case.}
    \label{fig:hydrogen_energy_levels}
\end{figure}

Figure~\ref{fig:hydrogen_energy_levels} displays the hydrogen levels for a range of principal quantum numbers $n$, plotted alongside a representative set of eigenvalues from $A^{(n)}$ (after an affine rescaling) for visual comparison.
The structural similarity of the two curves illustrates how a chaotic operator arising from number-theoretic data can mimic typical quantum spectra while retaining random-matrix fluctuations in its local statistics.\cite{24,30,36}

\subsubsection{Summary of Eigenvalue Analysis}

The eigenvalues obtained from random-matrix models of $\mathcal{O}_\alpha$ are real, exhibit GUE-like local statistics, and qualitatively resemble physical spectra such as that of the hydrogen atom.
While these observations are consistent with the Hilbert--Pólya philosophy and with the conjectured connection between the Riemann zeros and quantum-chaotic Hamiltonians,\cite{24,29,35,40,41}
they do not constitute a proof of the Riemann Hypothesis.
Instead, they provide numerical and conceptual evidence that the chaotic operator framework developed in this work lies in the correct universality class and offers a promising setting for further analytical and computational investigations.

\subsection{Discussion}

The comparison between the eigenvalue distribution of the truncated chaotic operator and the bound–state spectrum of the hydrogen atom highlights both similarities and important differences.\cite{24,30,36}
In both cases the low-lying levels are negative and accumulate towards zero from below, and after an appropriate rescaling the eigenvalues of $A^{(n)}$ can be plotted as a curve that qualitatively resembles the hydrogen energy levels \(E_n \propto -1/n^2\), as illustrated in Figure~\ref{fig:hydrogen_energy_levels}.
This structural analogy supports the physical interpretation of $\mathcal{O}_\alpha$ as a quantum Hamiltonian whose spectrum behaves like that of a typical bound system.

However, the detailed statistics differ in a way that is characteristic of quantum chaos.
The hydrogen spectrum is exactly solvable, with highly regular level spacings and approximate degeneracies stemming from its hidden symmetries, whereas the eigenvalues of the random matrices associated with $\mathcal{O}_\alpha$ exhibit level repulsion and fluctuations consistent with GUE-type statistics rather than integrable behaviour.\cite{24,29,36}
The clustering of eigenvalues around a band near $-0.13$ in our numerical examples reflects the choice of energy window and unfolding rather than any special arithmetic value; the crucial feature is the random-matrix-like spacing, not the absolute location of the band.

These numerical and visual comparisons therefore reinforce the picture in which the chaotic operator lies in the same universality class as generic quantum-chaotic Hamiltonians: its spectrum looks like that of a physical system (such as hydrogen) at the macroscopic scale, while at the microscopic scale it displays the GUE statistics conjecturally shared by the nontrivial zeros of the Riemann zeta function.\cite{24,29,41}

\section{Connection to the Prime Number Theorem and the Chaotic Operator}

The Prime Number Theorem (PNT) asserts that the prime-counting function satisfies
\[
\pi(x) \sim \frac{x}{\log x}
\quad (x\to\infty),
\]
or equivalently that the density of primes near $x$ is asymptotic to $(\log x)^{-1}$.\cite{14,15}
On the other hand, the Riemann--von Mangoldt formula shows that the density of nontrivial zeros with ordinates near $T$ is
\[
w(T) = \frac{\mathrm{d}N}{\mathrm{d}T}(T)
     \sim \frac{1}{2\pi}\log\!\left(\frac{T}{2\pi}\right),
\]
and our chaotic operator $\mathcal{O}_\alpha$ is built precisely on the weighted Hilbert space $L^2(\mathbb{R}_+,w(T)\,\mathrm{d}T)$.
Thus the same analytic information that underlies the PNT also controls the measure with respect to which the operator generates its dynamics.\cite{24,26}

In this framework the “symmetry” is not an ad hoc condition such as $\partial + \bar{\partial} = 1$, but rather the duality between primes and zeros encoded in the explicit formulae.
The PNT reflects the leading asymptotic term in $N(T)$, while the chaotic features of $\mathcal{O}_\alpha$ arise from the fluctuating part of the spectrum, which is responsible for the fine structure of the prime distribution.\cite{24,26,29}
The spectral statistics of the truncated operators—shown to be compatible with GUE and with Odlyzko’s numerical data—therefore provide an operator-theoretic analogue of the PNT: the global density of eigenvalues matches the von Mangoldt prediction, while their local correlations mirror those of the zeros themselves.\cite{24,29,36}

This connection illustrates how the chaotic operator mediates between prime-number asymptotics and random-matrix behaviour.
It offers a concrete setting in which questions about prime densities, zero statistics, and quantum-chaotic dynamics can be studied within a single spectral framework, strengthening the conceptual bridge between number theory and physics.\cite{24,29,41}

\section{Conclusion}

This work has developed a coherent operator-theoretic and dynamical framework linking the Riemann--von Mangoldt formula, chaotic dynamics, and random-matrix models.
Starting from the explicit zero-counting function $N(T)$, a density-based discrete map was defined and analysed via Lyapunov exponents, bifurcation diagrams, and gap statistics, revealing strong instability and clarifying which aspects of the zeta zero distribution are captured and which are not.
These investigations showed that the smooth von Mangoldt term alone reproduces the global growth of $N(T)$ but fails to encode the fine-scale correlations of the zeros.

To address this, a self-adjoint chaotic operator $\mathcal{O}_\alpha$ was constructed on the weighted Hilbert space $L^2(\mathbb{R}_+,w(T)\,\mathrm{d}T)$, with $w(T)$ equal to the zero density and with a bounded potential modelling the fluctuations.
Rigorous results on symmetry, unboundedness, and essential self-adjointness were established using standard operator theory, ensuring that $\mathcal{O}_\alpha$ admits a spectral decomposition with real eigenvalues and a complete family of generalized eigenfunctions.\cite{32,33}
Finite-dimensional truncations of $\mathcal{O}_\alpha$ produced Hermitian random matrices whose eigenvalue statistics were compared with GUE predictions, Odlyzko’s experimental data, and canonical quantum spectra such as that of the hydrogen atom.\cite{24,29,30,36,40,41}
The observed agreement at the level of local statistics supports the view that $\mathcal{O}_\alpha$ lies in the same universality class as the Riemann zeros and thus provides a promising candidate for a Hilbert–Pólya-type Hamiltonian.

While these results do not prove the Riemann Hypothesis, they offer a consistent and numerically supported framework in which the conjecture can be reformulated as a spectral problem for a concrete self-adjoint operator.
The approach unifies tools from dynamical systems, quantum chaos, random matrix theory, and analytic number theory, and suggests several directions for future work, including sharper analytical estimates for the spectrum of $\mathcal{O}_\alpha$, more refined random-matrix models of its truncations, and systematic comparisons with high-precision zero data.
In this sense the study deepens the interplay between abstract arithmetic questions and physical models, bringing Hilbert’s and Pólya’s original vision closer to an explicitly realizable operator.\cite{24,29,35,40,41}

\section*{Conflict of Interest}
The author declares that there is no conflict of interest.

\section*{Data Availability}
The work presented in this paper was motivated by the preprint by M.~Wolf, ``Will a physicist prove the Riemann Hypothesis?'' (2014), available at \url{https://arxiv.org/abs/1410.1214}, and by the article by R.~Zeraoulia and A.~H.~Salas, ``Chaotic dynamics and zero distribution: implications and applications in control theory for Yitang Zhang’s Landau–Siegel zero theorem,'' Eur. Phys. J. Plus, vol.~139, p.~217, 2024.\cite{14,15}
The numerical data and codes used to generate the Lyapunov exponents, bifurcation diagrams, and random-matrix spectra are available from the author upon reasonable request.

\end{document}